\documentclass{amsart}
\usepackage{url}
\newtheorem{thm}{Theorem}[section]
\newtheorem{cor}[thm]{Corollary}
\newtheorem{lem}[thm]{Lemma}

\theoremstyle{definition}

\theoremstyle{remark}

\begin{document}
\title[Slant submanifolds of quaternionic
space forms]{Slant submanifolds of quaternionic space forms}

\author[G.E. V\^{\i}lcu]{Gabriel Eduard V\^{\i}lcu}

\date{}
\maketitle

\abstract In this paper we establish some inequalities concerning
the $k$-Ricci curvature of a slant submanifold in a quaternionic
space form. We also obtain obstructions to the existence of
quaternionic slant immersions in quaternionic space forms
with unfull first normal bundle.\\ \\
{\bf Key Words:} Chen's invariant, scalar curvature, squared mean
curvature, $k$-Ricci curvature, quaternionic space form, slant
submanifold.\\ \\
{\bf AMS Mathematics Subject Classification:} 53C15, 53C25, 53C40.
\endabstract

\section{Introduction}

According to B.-Y. Chen \cite{CH5}, one of the most important
problems in submanifold theory is "\emph{to find simple
relationships between the main extrinsic invariants and the main
intrinsic invariants of a submanifold}". In \cite{CH4}, B.-Y. Chen
established a sharp inequality between the $k$-Ricci curvature, one
of the main intrinsic invariants, and the squared mean curvature,
the main extrinsic invariant, for a submanifold in a real space form
with arbitrary codimension. Also, in the same spirit, B.-Y. Chen
obtained an optimal inequality between the $k$-Ricci curvature and
the shape operator for submanifolds in real space forms. These
inequalities were further extended to many classes of submanifolds
in different ambient spaces: complex space forms \cite{MMO1,MMO2},
cosymplectic space forms \cite{LIU4,LIU3,YO2}, Sasakian space forms
\cite{CIO,KLY,MIH,TR}, locally conformal K\"{a}hler space forms
\cite{CAR,HON2}, generalized complex space forms \cite{HON,KIM,MA},
locally conformal almost cosymplectic manifolds \cite{KTC,YO1},
$(\kappa,\mu)$-contact space forms \cite{TRI}, Kenmotsu space forms
\cite{ARS,LIU}, $S$-space forms \cite{FER,KDT}.

In quaternionic setting, such inequalities were obtained for
quaternionic and totally-real submanifolds \cite{LIU1,LIU2,YO}. But
there are two classes of submanifolds which generalize both
quaternionic and totally real submanifolds of quaternionic
K\"{a}hler manifolds, with no inclusion between them: quaternionic
CR-submanifolds (see \cite{BCU}) and slant submanifolds (see
\cite{SH}). Some recent results concerning quaternionic
CR-submanifolds can be found in \cite{IIV,VILC} and an inequality
involving Ricci curvature and squared mean curvature for
quaternionic CR-submanifolds in quaternionic space forms was proved
in \cite{MSS}. On the other hand, some optimal inequalities
involving scalar curvature, Ricci curvature and squared mean
curvature for slant submanifolds in quaternionic space forms were
obtained recently in \cite{MHS,VIL}. The main purpose of this paper
is to obtain two kinds of inequalities for slant submanifolds in
quaternionic space forms: between the $k$-Ricci curvature and the
squared mean curvature and between the $k$-Ricci curvature and the
shape operator. Moreover, we investigate the existence of
quaternionic slant immersions in quaternionic space forms with
unfull first normal bundle.

\section{Preliminaries}

Let $\overline{M}$ be a differentiable manifold and assume that
there is a rank 3-subbundle $\sigma$ of $End(T\overline{M})$ such
that a local basis $\lbrace{J_1,J_2,J_3}\rbrace$ exists on sections
of $\sigma$ satisfying for all $\alpha\in\{1,2,3\}$:
       \begin{equation}\label{1}
       J_\alpha^2=-Id,
       J_{\alpha}J_{\alpha+1}=-J_{\alpha+1}J_{\alpha}=J_{\alpha+2},
       \end{equation}
where the indices are taken from $\{1,2,3\}$ modulo 3. Then the
bundle $\sigma$ is called an almost quaternionic structure on $M$
and $\lbrace{J_1,J_2,J_3}\rbrace$ is called a canonical local basis
of $\sigma$. Moreover, $(\overline{M},\sigma)$ is said to be an
almost quaternionic manifold. It is easy to see that any almost
quaternionic manifold is of dimension $4m$.

A Riemannian metric $\overline{g}$ on $\overline{M}$ is said to be
adapted to the almost quaternionic structure $\sigma$ if it
satisfies:
 \begin{equation}\label{2}
         \overline{g}(J_\alpha X,J_\alpha Y)=
         \overline{g}(X,Y),\forall\alpha\in\{1,2,3\}
         \end{equation}
for all vector fields $X$,$Y$ on $\overline{M}$ and any canonical
local basis $\lbrace{J_1,J_2,J_3}\rbrace$ of $\sigma$. Moreover,
$(\overline{M},\sigma,\overline{g})$ is said to be an almost
quaternionic Hermitian manifold.

If the bundle $\sigma$ is parallel with respect to the Levi-Civita
connection $\overline{\nabla}$ of $\overline{g}$, then
$(\overline{M},\sigma,\overline{g})$ is said to be a quaternionic
K\"{a}hler manifold. Equivalently, locally defined 1-forms
$\omega_1,\omega_2,\omega_3$ exist such that we have for all
$\alpha\in\{1,2,3\}$:
    \begin{equation}\label{3}
     \overline{\nabla}_XJ_{\alpha}=\omega_{\alpha+2}(X)J_{\alpha+1}-
    \omega_{\alpha+1}(X)J_{\alpha+2},
    \end{equation}
for any vector field $X$ on $\overline{M}$, where the indices are
taken from $\{1,2,3\}$ modulo 3 (see \cite{ISH}).

Let $(\overline{M},\sigma,\overline{g})$ be a quaternionic
K\"{a}hler manifold and let $X$ be a non-null vector on
$\overline{M}$. Then the 4-plane spanned by $\{X,J_1X,J_2X,J_3X\}$,
denoted by $Q(X)$, is called a quaternionic 4-plane. Any 2-plane in
$Q(X)$ is called a quaternionic plane. The sectional curvature of a
quaternionic plane is called a quaternionic sectional curvature. A
quaternionic K\"{a}hler manifold is a quaternionic space form if its
quaternionic sectional curvatures are equal to a constant, say $c$.
It is well-known that a quaternionic K\"{a}hler manifold
$(\overline{M},\sigma,\overline{g})$ is a quaternionic space form,
denoted $\overline{M}(c)$, if and only if its curvature tensor is
given by (see \cite{ISH}):
     \begin{eqnarray}\label{4}
       \overline{R}(X,Y)Z&=&\frac{c}{4}\lbrace\overline{g}(Z,Y)X-
       \overline{g}(X,Z)Y+\sum\limits_{\alpha=1}^3
        [\overline{g}(Z,J_\alpha Y)J_\alpha X-\nonumber\\
       &&-\overline{g}(Z,J_\alpha X)J_\alpha Y+
        2\overline{g}(X,J_\alpha Y)J_\alpha Z]\rbrace
      \end{eqnarray}
for all vector fields $X,Y,Z$ on $\overline{M}$ and any local basis
$\lbrace{J_1,J_2,J_3}\rbrace$ of $\sigma.$

For a submanifold $M$ of a quaternion K\"{a}hler manifold
$(\overline{M},\sigma,\overline{g})$, we denote by $g$ the metric
tensor induced on $M$. If $\nabla$ is the covariant differentiation
induced on $M$, the Gauss and Weingarten formulas are given by:
\begin{equation}\label{5}
       \overline{\nabla}_XY=\nabla_XY+h(X,Y), \forall X,Y \in
\Gamma(TM)
       \end{equation}
and
\begin{equation}\label{6}
       \overline{\nabla}_XN=-A_NX+\nabla_{X}^{\perp}N, \forall X\in
\Gamma(TM), \forall N\in \Gamma(TM^\perp)
       \end{equation}
where $h$ is the second fundamental form of $M$, $\nabla^\perp$ is
the connection on the normal bundle and $A_N$ is the shape operator
of $M$ with respect to $N$. The shape operator $A_N$ is related to
$h$ by:
\begin{equation}\label{7}
       g(A_NX,Y)=\overline{g}(h(X,Y),N),
       \end{equation}
for all $ X,Y\in \Gamma(TM)$ and $N\in \Gamma(TM^\perp)$.

If we denote by $\overline{R}$ and $R$ the curvature tensor fields
of $\overline{\nabla}$ and $\nabla$ we have the Gauss equation:
\begin{eqnarray}\label{8}
&&\overline{R}(X,Y,Z,W)=R(X,Y,Z,W)+\overline{g}(h(X,W),h(Y,Z))-\overline{g}(h(X,Z),h(Y,W)),
       \end{eqnarray}
for all $X,Y,Z,W\in \Gamma(TM)$.

For the second fundamental form $h$, we define the covariant
derivative $\overline{\nabla}h$ of $h$ with respect to the
connection on $TM\oplus T^\perp M$ by
\begin{equation}\label{8.1}
(\overline{\nabla}_Xh)(Y,Z)=D_X(h(Y,Z))-h(\nabla_XY,Z)-h(Y,\nabla_XZ),
\end{equation}
where $D$ is the linear connection induced on the normal bundle of
$M$ in $\overline{M}$. Then the equation of Codazzi is given by
\begin{equation}\label{8.2}
(\overline{R}(X,Y)Z)^\perp=(\overline{\nabla}_Xh)(Y,Z)-(\overline{\nabla}_Yh)(X,Z).
\end{equation}

If $\{e_1,...,e_n\}$ is an orthonormal basis of $T_pM$ and
$\{e_{n+1},...,e_{4m}\}$ is an orthonormal basis of $T_p^\perp M$,
where $p\in M$, we denote by $H$ the mean curvature vector, that is
$$H(p)=\frac{1}{n}\sum_{i=1}^{n}h(e_i,e_i).$$

Also, we set
$$h_{ij}^r=g(h(e_i,e_j),e_r),\ i,j\in\{1,...,n\},\ r\in\{n+1,...,4m\}$$
and
$$||h||^2(p)=\sum_{i,j=1}^{n}g(h(e_i,e_j),h(e_i,e_j)).$$

A submanifold $M$ of a quaternionic K\"{a}hler manifold
$\overline{M}$ is called a quaternionic submanifold (resp. totally
real submanifold) if each tangent space of $M$ is carried into
itself (resp. into the normal space) by each section in $\sigma$.
Recently, \c Sahin \cite{SH} introduced the slant submanifolds of
quaternionic K\"{a}hler manifolds, as a natural generalization of
both quaternionic and totally real submanifolds. A submanifold $M$
of a quaternionic K\"{a}hler manifold $\overline{M}$ is said to be a
slant submanifold if for each non-null vector $X$ tangent to $M$ at
$p$, the angle $\theta(X)$ between $J_\alpha(X)$ and $T_pM$,
$\alpha\in\{1,2,3\}$, is constant, i.e. it does not depend on choice
of $p\in M$ and $X\in T_pM$. We can easily see that quaternionic
submanifolds are slant submanifolds with $\theta=0$ and totally-real
submanifolds are slant submanifolds with $\theta=\frac{\pi}{2}$. A
slant submanifold of a quaternionic K\"{a}ler manifold is said to be
proper (or $\theta$-slant proper) if it is neither quaternionic nor
totally real.

If $M$ is a slant submanifold of a quaternionic K\"{a}hler manifold
$\overline{M}$, then for any $X\in\Gamma(TM)$ we have the
decomposition
\begin{equation}\label{D1}
J_\alpha X=P_{\alpha}X+F_{\alpha}X,
\end{equation}
where $P_{\alpha}X$ denotes the tangential component of
$J_{\alpha}X$ and $F_{\alpha}X$ denotes the normal component of
$J_{\alpha}X$.

Similarly for any $U\in\Gamma(TM^\perp)$ we have
\begin{equation}\label{D2}
J_\alpha U=B_{\alpha}U+C_{\alpha}U,
\end{equation}
where $B_{\alpha}U$ is the tangential component of $J_{\alpha}U$ and
$C_{\alpha}U$ is the normal component of $J_{\alpha}U$.

We recall now the following results which we shall need in the
sequel.

\begin{thm} \cite{SH}
Let $M$ be a submanifold of a quaternionic K\"{a}hler manifold
$\overline{M}$. Then $M$ is slant if and only if there exists a
constant $\lambda\in[-1,0]$ such that:
\begin{equation}\label{9}
P_{\beta}P_{\alpha}X=\lambda X,\ \ \forall X\in\Gamma(TM),\
\alpha,\beta\in\{1,2,3\}.
\end{equation}
Furthermore, in such case, if $\theta$ is the slant angle of $M$,
then it satisfies $\lambda=-{\rm cos}^2\theta$.
\end{thm}

\begin{cor} \cite{SH}
Let $M$ be a slant submanifold of a quaternionic K\"{a}hler manifold
$\overline{M}$, with slant angle $\theta$. Then we have
\begin{equation}\label{a.a}
P^2_\alpha X=-{\rm cos}^2\theta X,
\end{equation}
\begin{equation}\label{a.b}
B_\alpha F_\alpha X=-{\rm sin}^2\theta X,
\end{equation}
for any $X\in\Gamma(TM)$ and $\alpha\in\{1,2,3\}$.
\end{cor}

From the above Theorem we deduce that if $M$ is a $\theta$-slant
submanifold of a quaternionic K\"{a}hler manifold $\overline{M}$,
then we have for any $X,Y\in \Gamma(TM)$:
\begin{equation}\label{10}
g(P_\alpha X,P_\beta Y)={\rm cos}^2\theta g(X,Y),\
\alpha,\beta\in\{1,2,3\}
\end{equation}
and
\begin{equation}\label{10b}
\overline{g}(F_\alpha X,F_\beta Y)={\rm sin}^2\theta g(X,Y),\
\alpha,\beta\in\{1,2,3\}.
\end{equation}

Moreover, we can remark that every proper slant submanifold of a
quaternionic K\"{a}hler manifold is of even dimension $n=2s$,
because we can choose a canonical orthonormal local frame
$\{e_1,{\rm sec}\theta P_\alpha e_1,...,e_s,{\rm sec}\theta P_\alpha
e_s\}$ of $T_pM$, $p\in M$, called an adapted slant frame, where
$\alpha$ is settled in $\{1,2,3\}$.

For an $n$-dimensional Riemanian manifold $(M,g)$ we denote by
$K(\pi)$ the sectional curvature of $M$ associated with a plane
section $\pi\subset T_pM,\ p\in M$. If $\{e_1,...,e_n\}$ is an
orthonormal basis of the tangent space $T_pM$, the scalar curvature
$\tau$ at $p$ is defined by
\begin{equation}\label{11}
\tau(p)=\sum_{1\leq i<j\leq n}K_{ij},
\end{equation}
where $K_{ij}$ denotes the sectional curvatures of the 2-plane
section spanned by $e_i$ and $e_j$.

For a $k$-plane section $L$ of $T_pM$, $p\in M$, and $X$ a unit
vector in $L$, we choose an orthonormal basis $\{e_1,...,e_k\}$ of
$L$ such that $e_1=X$. The Ricci curvature of $L$ at $X$, denoted
$Ric_L(X)$, is defined by
\begin{equation}\label{12}
Ric_L(X)=\sum_{j=2}^{k}K_{1j}.
\end{equation}
We note that such a curvature is called a $k$-Ricci curvature. The
scalar curvature of a $k$-plane section $L$ is given by
\begin{equation}\label{11bis}
\tau(L)=\sum_{1\leq i<j\leq k}K_{ij},
\end{equation}

For an integer $k$, $2\leq k\leq n$, B.-Y. Chen introduced a
Riemannian invariant $\Theta_k$ defined by
\begin{equation}\label{13}
\Theta_k(p)=\frac{1}{k-1}inf\{Ric_L(X)|L,X\},\ p\in M,
\end{equation}
where $L$ runs over all $k$-plane sections in $T_pM$ and $X$ runs
over all unit vectors in $L$ (see e.g. \cite{CH01}).

\section{$k$-Ricci curvature and the squared mean curvature}

\begin{thm}\label{1.1}
Let $M^n$ be a $\theta$-slant proper submanifold of a quaternionic
space form $\overline{M}^{4m}(c)$. Then, for any $p\in M$ and any
integer $k$, $2\leq k\leq n$, one has:
\begin{equation}\label{14}
||H||^2(p)\geq \Theta_k(p)-\frac{c}{4}\left(1+\frac{9}{n-1}{\rm
cos}^2\theta\right).
\end{equation}
\end{thm}
\begin{proof}
We choose an adapted slant basis of $T_pM$ at $p\in M$:
$$\{e_1,e_2={\rm sec}\theta P_\alpha e_1,...,e_{2s-1},e_{2s}={\rm sec}\theta
P_\alpha e_{2s-1}\},$$  where $2s=n$, and $\{e_{n+1},...,e_{4m}\}$
an orthonormal basis of $T_p^\perp M$, such that the normal vector
$e_{n+1}$ is in the direction of the mean curvature vector $H(p)$
and $\{e_1,...,e_n\}$ diagonalize the shape operator $A_{n+1}$.

Taking now $X=Z=e_i$, $Y=W=e_j$ in the equation of Gauss (\ref{8}),
by summing and using (\ref{4}), we obtain:
\begin{equation}\label{15}
n^2||H||^2(p)=2\tau(p)+||h||^2(p)-\frac{n(n-1)c}{4}-\frac{3c}{4}\sum_{\beta=1}^3\sum_{i,j=1}^{n}g^2(P_\beta
e_i,e_j).
\end{equation}

On the other hand, because $\{e_1,...,e_{2s}\}$ is an adapted slant
basis of $T_pM$, using (\ref{9}) and (\ref{10}) we can see that we
have:
\begin{equation}\label{16}
g^2(P_\beta e_i,e_{i+1})=g^2(P_\beta e_{i+1},e_i)={\rm
cos}^2\theta,\ {\rm for}\ i=1,3,...,2s-1
\end{equation}
and
\begin{equation}\label{17}
g(P_\beta e_i,e_j)=0,\ {\rm for}\
(i,j)\not\in\{(2l-1,2l),(2l,2l-1)|l\in\{1,2,....,s\}\}.
\end{equation}

From (\ref{15}), (\ref{16}) and (\ref{17}) we derive:
\begin{equation}\label{18}
n^2||H||^2(p)=2\tau(p)+||h||^2(p)-\frac{c}{4}\left[n(n-1)+9n{\rm
cos}^2\theta\right].
\end{equation}

On the other hand, due to the choosing of the basis of $T_pM$ and
$T_p^\perp M$, the shape operators have the following forms:
\begin{equation}\label{19}
A_{n+1}= \left( \begin{array}{ccccc}
a_1 & 0 & 0 & ... & 0\\
0 & a_2 & 0 & ... & 0\\
0 & 0 & a_3 & ... & 0\\
\vdots & \vdots & \vdots & \ddots & \vdots\\
0 & 0 & 0 & ... &  a_n \end{array} \right),
\end{equation}
\begin{equation}\label{20}
A_r=(h_{ij}^r)_{i,j=\overline{1,n}},\
traceA_r=\sum_{i=1}^{n}h_{ii}^r=0,\ \forall r\in\{n+2,...,4m\}.
\end{equation}

Now, using (\ref{19}) and (\ref{20}) in (\ref{18}) we obtain:
\begin{equation}\label{21}
n^2||H||^2(p)=2\tau(p)+\sum_{i=1}^{n}a_i^2+\sum_{r=n+2}^{4m}\sum_{i,j=1}^{n}(h_{ij}^{r})^2-\frac{c}{4}\left[n(n-1)+9n{\rm
cos}^2\theta\right].
\end{equation}

On the other hand, because we have the inequality
$$
(n-1)\sum_{i=1}^{n}a_i^2\geq2\sum_{i<j}a_ia_j,
$$
from
$$
n^2||H||^2(p)=\left(\sum_{i=1}^na_i\right)^2=\sum_{i=1}^{n}a_i^2+2\sum_{1\leq
i<j\leq n}a_ia_j
$$
we derive
\begin{equation}\label{22}
\sum_{i=1}^{n}a_i^2\geq n||H||^2(p).
\end{equation}

Using now (\ref{22}) in (\ref{21}) we obtain:
\begin{equation}\label{23}
n(n-1)||H||^2(p)\geq2\tau(p)-\frac{c}{4}\left[n(n-1)+9n{\rm
cos}^2\theta\right].
\end{equation}

But, from (\ref{11}) and (\ref{11bis}), it follows that for any
$k$-plane section $L_{i_1...i_k}$ spanned by
$\{e_{i_1},...,e_{i_k}\}$, one has:
\begin{equation}\label{24}
\tau(L_{i_1...i_k})=\frac{1}{2}\sum_{i\in\{i_1,...,i_k\}}Ric_{L_{i_1...i_k}}(e_i)
\end{equation}
and
\begin{equation}\label{25}
\tau(p)=\frac{(k-2)!(n-k)!}{(n-2)!}\sum_{1\leq i_1<...<i_k\leq
n}\tau(L_{i_1...i_k}).
\end{equation}

From (\ref{24}) and (\ref{25}) we obtain:
\begin{equation}\label{26}
\tau(p)\geq \frac{n(n-1)}{2}\cdot\Theta_k(p)
\end{equation}
and finally, from (\ref{23}) and (\ref{26}) one derives (\ref{14}).
\end{proof}

Applying Theorem \ref{1.1} we may obtain, as a particular case, the
corresponding inequality for totally-real submanifolds in
quaternionic space forms, established in \cite{LIU2}.

\begin{cor}
Let $M^n$ be a totally-real submanifold of a quaternionic space form
$\overline{M}^{4m}(c)$. Then, for any $p\in M$ and any integer $k$,
$2\leq k\leq n$, one has:
\begin{equation}
||H||^2(p)\geq \Theta_k(p)-\frac{c}{4}.
\end{equation}
\end{cor}

\section{$k$-Ricci curvature and shape operator}

\begin{thm}\label{4.1}
Let $x:M\rightarrow\overline{M}^{4m}(c)$ be an isometric immersion
of an $n$-dimensional $\theta$-slant proper submanifold $M$ into a
$4m$-dimensional quaternionic space form $\overline{M}(c)$. Then,
for any $p\in M$ and any integer $k$, $2\leq k\leq n$, one has:

i. If $\Theta_k(p)\neq\frac{c}{4}\left(1+\frac{9}{n-1}{\rm
cos}^2\theta\right)$, then the shape operator at the mean curvature
satisfies
\begin{equation}
A_H>
\frac{n-1}{n}\left[\Theta_k(p)-\frac{c}{4}\left(1+\frac{9}{n-1}{\rm
cos}^2\theta\right)\right]I_n,
\end{equation}
at $p$, where $I_n$ denotes the identity map of $T_pM$.

ii. If $\Theta_k(p)=\frac{c}{4}\left(1+\frac{9}{n-1}{\rm
cos}^2\theta\right)$, then $A_H\geq 0$ at $p$.

iii. A unit vector $X\in T_pM$ satisfies
\begin{equation}\label{28}
A_HX=
\frac{n-1}{n}\left[\Theta_k(p)-\frac{c}{4}\left(1+\frac{9}{n-1}{\rm
cos}^2\theta\right)\right]X
\end{equation}
if and only if $\Theta_k(p)=\frac{c}{4}\left(1+\frac{9}{n-1}{\rm
cos}^2\theta\right)$ and $X$ belongs to the relative null space of
$M$ at $p$:
$$\mathcal{N}_p=\{Z\in T_pM|h(Z,Y)=0, \forall Y\in T_pM\}.$$

iv. The identity
\begin{equation}\label{29}
A_H=
\frac{n-1}{n}\left[\Theta_k(p)-\frac{c}{4}\left(1+\frac{9}{n-1}{\rm
cos}^2\theta\right)\right]I_n,
\end{equation}
holds at $p$ if and only if $p$ is a totally geodesic point.
\end{thm}
\begin{proof}
i. We choose an adapted slant basis of $T_pM$ at $p\in M$:
$$\{e_1,e_2={\rm sec}\theta P_\alpha e_1,...,e_{2s-1},e_{2s}={\rm sec}\theta
P_\alpha e_{2s-1}\},$$  where $2s=n$, and $\{e_{n+1},...,e_{4m}\}$
an orthonormal basis of $T_p^\perp M$, such that the normal vector
$e_{n+1}$ is in the direction of the mean curvature vector $H(p)$
and $\{e_1,...,e_n\}$ diagonalize the shape operator $A_{n+1}$.
Consequently, the shape operators have the forms (\ref{19}) and
(\ref{20}).

One can distinguishes two cases:

Case I: $H(p)=0$. In this situation it follows from (\ref{14}) that
$\Theta_k(p)\neq\frac{c}{4}\left(1+\frac{9}{n-1}{\rm
cos}^2\theta\right)$ and the conclusion follows.

Case II: $H(p)\neq0$. Taking $X=Z=e_i$ and $Y=W=e_j$ in the Gauss
equation and using (\ref{4}), we obtain:
\begin{eqnarray}\label{30}
a_ia_j=K_{ij}-\frac{c}{4}\left[1+3\sum_{\beta=1}^{3}g^2(P_\beta
e_i,e_j)\right]-\sum_{r=n+2}^{4m}\left[h_{ii}^rh_{jj}^r-(h_{ij}^r)^2\right].
\end{eqnarray}

From (\ref{30}) we derive:
\begin{eqnarray}\label{31}
a_1(a_{i_2}+...+a_{i_k})&=&Ric_{L_{1i_2...i_k}}(e_1)-\frac{(k-1)c}{4}-\frac{3c}{4}\sum_{\beta=1}^{3}\sum_{j=2}^{k}g^2(P_\beta
e_1,e_{i_j})\nonumber\\&&-\sum_{r=n+2}^{4m}\sum_{j=2}^{k}\left[h_{11}^rh_{i_ji_j}^r-(h_{1i_j}^r)^2\right]
\end{eqnarray}
which implies
\begin{eqnarray}\label{32}
a_1(a_2+...+a_n)&=&\frac{(k-2)!(n-k)!}{(n-2)!}\sum_{2\leq
i_2<...<i_k\leq
n}Ric_{L_{1i_2...i_k}}(e_1)-\frac{(n-1)c}{4}\nonumber\\&&-\frac{3c}{4}\sum_{\beta=1}^{3}\sum_{j=2}^{n}g^2(P_\beta
e_1,e_j)+\sum_{r=n+2}^{4m}\sum_{j=1}^{n}(h_{1j}^r)^2.
\end{eqnarray}
and taking into account (\ref{13}), we obtain:
\begin{eqnarray}\label{33}
a_1(a_2+...+a_n)\geq(n-1)\theta_k(p)-\frac{(n-1)c}{4}-\frac{3c}{4}\sum_{\beta=1}^{3}\sum_{j=2}^{n}g^2(P_\beta
e_1,e_j).
\end{eqnarray}

Using (\ref{16}) and (\ref{17}) in (\ref{33}) we obtain:
\begin{eqnarray}\label{34}
a_1(a_2+...+a_n)\geq(n-1)\Theta_k(p)-\frac{(n-1)c}{4}-\frac{9c}{4}{\rm
cos}^2\theta
\end{eqnarray}
and we find:
\begin{eqnarray}\label{35}
a_1(a_1+a_2+...+a_n)&=&a_1^2+a_1(a_2+...+a_n)\nonumber\\
&\geq&(n-1)\left[\Theta_k(p)-\frac{c}{4}\left(1+9{\rm
cos}^2\theta\right)\right].
\end{eqnarray}

Similar inequalities hold when the index 1 is replaced by
$j\in\{2,...,n\}$. Hence, we have
\begin{eqnarray}\label{36}
a_j(a_1+a_2+...+a_n)\geq
(n-1)\left[\Theta_k(p)-\frac{c}{4}\left(1+9{\rm
cos}^2\theta\right)\right],
\end{eqnarray}
for all $j\in\{1,...,n\}$, and because $n||H||=a_1+...+a_n$ we find
\begin{equation}\label{37}
A_H\geq
\frac{n-1}{n}\left[\Theta_k(p)-\frac{c}{4}\left(1+\frac{9}{n-1}{\rm
cos}^2\theta\right)\right]I_n.
\end{equation}

We remark that the equality does not hold because we are in the case
$H(p)\neq0$.

ii. The statement is clear from i.

iii. If $X\in T_pM$ is a unit vector such that (\ref{28}) holds,
then we have equalities both in (\ref{33}) and (\ref{35}).
Consequently, we obtain $a_1=0$ and $h_{1j}^r=0$, for all
$j\in\{1,...,n\}$ and $r\in\{n+2,...,4m\}$, which implies
$\Theta_k(p)=\frac{c}{4}\left(1+\frac{9}{n-1}{\rm
cos}^2\theta\right)$ and $X\in\mathcal{N}_p$. The converse part is
clear.

iv. The equality (\ref{32}) holds for any $X\in T_pM$ if and only if
$\mathcal{N}_p=TpM$, i.e. $p$ is a totally geodesic point. This
completes the proof of the theorem.
\end{proof}

\begin{cor}
Let $x:M\rightarrow\overline{M}^{4m}(c)$ be an isometric immersion
of an $n$-dimensional totally-real submanifold $M$ into a
$4m$-dimensional quaternionic space form $\overline{M}(c)$. Then,
for any $p\in M$ and any integer $k$, $2\leq k\leq n$, one has:

i. If $\Theta_k(p)\neq\frac{c}{4}$, then the shape operator at the
mean curvature satisfies
\begin{equation}
A_H> \frac{n-1}{n}\left[\Theta_k(p)-\frac{c}{4}\right]I_n,
\end{equation}
at $p$, where $I_n$ denotes the identity map of $T_pM$.

ii. If $\Theta_k(p)=\frac{c}{4}$, then $A_H\geq 0$ at $p$.

iii. A unit vector $X\in T_pM$ satisfies
\begin{equation}
A_HX= \frac{n-1}{n}\left[\Theta_k(p)-\frac{c}{4}\right]X
\end{equation}
if and only if $\Theta_k(p)=\frac{c}{4}$ and $X\in\mathcal{N}_p$.

iv. The identity
\begin{equation}
A_H= \frac{n-1}{n}\left[\Theta_k(p)-\frac{c}{4}\right]I_n,
\end{equation}
holds at $p$ if and only if $p$ is a totally geodesic point.
\end{cor}

\section{Quaternionic slant submanifolds with unfull first normal bundle}

Let $M$ be a submanifold isometrically immersed in  a Riemannian
manifold $(\overline{M},\overline{g})$. If $p$ is a point of $M$,
then the first normal space at $p$ is defined to be ${\rm Im} h_p$,
the image space of the second fundamental form $h$ at $p$. Moreover,
${\rm Im} h$ is called the first normal bundle of $M$ in
$\overline{M}$. The submanifold is said to have full first normal
bundle if ${\rm Im} h_p=T_pM^\perp$, for any $p\in M$ (see
\cite{CH0,CH1}).

The existence of K\"{a}hlerian slant submanifolds of smallest
possible codimension in complex space forms, having unfull first
normal bundle, has been investigated in \cite{LIWU}. Next we'll
study this problem in the context of slant submanifolds in
quaternionic space forms. The quaternionic version of K\"{a}hlerian
slant submanifolds has been introduced in \cite{SH}, under the name
of quaternionic slant submanifolds. Therefore a proper slant
submanifold $M$ of a quaternionic K\"{a}hler manifold
$(\overline{M},\sigma,\overline{g})$ is said to be quaternionic
slant submanifold if it satisfies the condition
\begin{equation}
     \overline{\nabla}_XP_{\alpha}=\omega_{\alpha+2}(X)P_{\alpha+1}-
    \omega_{\alpha+1}(X)P_{\alpha+2},
\end{equation}
for any vector field $X$ on $\overline{M}$, where the indices are
taken from $\{1,2,3\}$ modulo 3.

We have the following characterization of quaternionic slant
submanifolds.

\begin{thm}\label{5.1} \cite{SH}
Let $M$ be a proper slant submanifold of a quaternionic K\"{a}hler
manifold $\overline{M}$. Then $M$ is quaternionic slant submanifold
if and only if
\begin{equation}\label{51}
A_{F\alpha Y}Z=A_{F\alpha Z}Y
\end{equation}
for all $Y,Z\in\Gamma(TM)$ and $\alpha\in\{1,2,3\}$.
\end{thm}

\begin{lem}
Let $M$ be a slant submanifold of a quaternionic K\"{a}hler manifold
$\overline{M}$. Then we have
\begin{equation}\label{L1}
P^2_\alpha=-Id-B_\alpha F_\alpha,
\end{equation}
\begin{equation}\label{L2}
C_\alpha F_\alpha+F_\alpha P_\alpha=0,
\end{equation}
\begin{equation}\label{L3}
C^2_\alpha=-Id-F_\alpha B_\alpha,
\end{equation}
\begin{equation}\label{L4}
P_\alpha B_\alpha+B_\alpha C_\alpha=0
\end{equation}
for $\alpha\in\{1,2,3\}$.
\end{lem}
\begin{proof}
For any $X\in\Gamma(TM)$, taking into account (\ref{1}) and
(\ref{D1}), we derive
\[
-X=J^2_\alpha X=P^2_\alpha X+F_\alpha P_\alpha X+B_\alpha F_\alpha
X+C_\alpha F_\alpha X.
\]
Equating the tangent and normal parts of both the sides we obtain
(\ref{L1}) and (\ref{L2}).

On the other hand, taking into account (\ref{1}) and (\ref{D2}) we
conclude that for any $U\in\Gamma(TM^\perp)$ we have
\[
-U=J^2_\alpha U=P_\alpha B_\alpha U+F_\alpha B_\alpha U+B_\alpha
C_\alpha U+C^2_\alpha U.
\]
Equating now the tangent and normal parts of both the sides we
obtain (\ref{L3}) and (\ref{L4}).
\end{proof}

\begin{lem}\label{5.3}
Let $M$ be a $\theta$-slant proper submanifold of a quaternionic
K\"{a}hler manifold $\overline{M}$. Then for any vectors $U,V\in
T_pM^\perp$, $p\in M$, we have
\begin{equation}\label{56}
\overline{g}(C_\alpha U,C_\alpha V)={\rm cos}^2\theta
\overline{g}(U,V),\
 \alpha=1,2,3.
\end{equation}
\end{lem}
\begin{proof}
Because $M$ is a $\theta$-slant proper submanifold of
$\overline{M}$, it follows that there exist $X_\alpha,Y_\alpha\in
T_pM$ such that $U=F_\alpha X_\alpha,\ V=F_\alpha Y_\alpha$. Then,
by using (\ref{10}), (\ref{10b}) and (\ref{L2}), we derive
\begin{eqnarray*}
\overline{g}(C_\alpha U,C_\alpha V)&=&\overline{g}(C_\alpha F_\alpha
X_\alpha,C_\alpha F_\alpha Y_\alpha)
=\overline{g}(F_\alpha P_\alpha X_\alpha,F_\alpha P_\alpha Y_\alpha)\\
&=&{\rm sin}^2\theta \overline{g}(P_\alpha X_\alpha,P_\alpha Y_\alpha)={\rm sin}^2\theta {\rm cos}^2\theta \overline{g}(X_\alpha,Y_\alpha)\\
&=&{\rm cos}^2\theta \overline{g}(F_\alpha X_\alpha,F_\alpha
Y_\alpha)={\rm cos}^2\theta \overline{g}(U,V).
\end{eqnarray*}
\end{proof}

From Theorem \ref{5.1} and Lemma \ref{5.3}, using the same
techniques as in \cite{LIWU}, we can state now the following result.
\begin{lem}\label{5.4}
Let $M$ be a quaternionic slant submanifold of a quaternionic
K\"{a}hler manifold $\overline{M}$. Then \[B_\alpha({\rm Im}
h_p)^\perp=\mathcal{N}_p,\ \alpha=1,2,3,\] where $(Im h_p)^\perp$
denotes the orthogonal complementary subspace of $Im h_p$ in
$T_pM^\perp$ and $\mathcal{N}_p$ is the relative null space of $M$
at $p$.
\end{lem}
\begin{proof}
For $Z\in B_\alpha(Im h_p)^\perp$ it follows that there exists $U\in
(Im h_p)^\perp$ such that $Z=B_\alpha U$. Then, by using (\ref{7}),
(\ref{51}), (\ref{L3}) and (\ref{56}), we obtain for all vector
$X,Y\in T_pM$ and $\alpha=1,2,3$:
\begin{eqnarray*}
\overline{g}(h(X,Z),F_\alpha Y)&=&g(A_{F_\alpha Y}Z,X)=g(A_{F_\alpha
Z}Y,X)\\
&=&\overline{g}(h(X,Y),F_\alpha Z)=\overline{g}(h(X,Y),F_\alpha
B_\alpha U)\\
&=&{\rm sin}^2\theta\overline{g}(h(X,Y),U)=0.
\end{eqnarray*}

Therefore it follows that $h(X,Z)=0$, for any $X\in T_pM$ and thus
we obtain $Z\in \mathcal{N}_p$.

If we take now $Z\in \mathcal{N}_p$, it is clear that for any
$X,Y\in T_pM$ and $\alpha=1,2,3$ we have
\[
\overline{g}(h(X,Y),F_\alpha Z)=\overline{g}(h(Z,X),F_\alpha Y)=0.
\]

Thus it follows $F_\alpha Z\in ({\rm Im} h_p)^\perp$ and therefore
we derive
\begin{equation}\label{57}
B_\alpha F_\alpha Z\in B_\alpha({\rm Im} h_p)^\perp.
\end{equation}

From (\ref{a.b}) and (\ref{57}) we conclude that $Z\in B_\alpha({\rm
Im} h_p)^\perp$ and the proof is now complete.
\end{proof}

\begin{thm}
Let $x:M\rightarrow\overline{M}(c)$ be an isometric immersion of a
quaternionic slant submanifold $M$ of minimal codimension into a
quaternionic space form $\overline{M}(c)$. If the first normal
bundle is not full, then $c=0$.
\end{thm}
\begin{proof}
First of all we remark that if the dimension of $\overline{M}(c)$ is
$4m$, then the minimal codimension of a proper slant submanifold $M$
of $\overline{M}(c)$ is $2m$; in this case we can choose an adapted
slant basis of $T_pM$ at $p\in M$:
\[\{e_1,e_2={\rm sec}\theta P_\alpha e_1,...,e_{2m-1},e_{2m}={\rm sec}\theta
P_\alpha e_{2m-1}\},\] and an orthonormal basis of $T_p^\perp M$:
\[\{e_{2m+1}={\rm cosec}\theta F_\alpha e_1,e_{2m+2}={\rm cosec}\theta F_\alpha e_2,...,e_{4m}={\rm cosec}\theta F_\alpha e_{2m}\},\]
where $\alpha$ is settled in $\{1,2,3\}$.

Moreover, if the first normal bundle is not full, then it follows
that there exists a unit normal vector $U\in T_pM^\perp$ at a point
$p\in M$ such that $\overline{g}(h(X,Y),U)=0$, for any vector
$X,Y\in T_pM$ and without loss of generality we can suppose
$e_{4m}=U$. Applying Lemma \ref{5.4} it follows $B_\alpha
e_{4m}\in\mathcal{N}_p$, $\alpha=1,2,3$, and from (\ref{a.b}) we
conclude $e_{2m}\in\mathcal{N}_p$. Thus we have
\begin{equation}\label{58}
h(e_i,e_{2m})=0,\ i=1,...,2m-1.
\end{equation}

By using now (\ref{8.1}) and (\ref{58}) in (\ref{8.2}) we obtain for
$i=1,...,2m-1$:
\[
(\overline{R}(e_i,e_{2m})e_{2m})^\perp=h(e_i,\nabla_{e_{2m}}e_{2m})
\]
and taking into account (\ref{58}) and the definition of the
Christoffel symbols $\Gamma_{ij}^k$:
\[
\nabla_{e_i}e_{j}=\sum_{k=1}^{2m}\Gamma_{ij}^ke_k
\]
we obtain
\begin{equation}\label{59}
(\overline{R}(e_i,e_{2m})e_{2m})^\perp=\sum_{1\leq
k,l<2m}\Gamma_{2m2m}^kh_{ik}^{2m+l}e_{2m+l}.
\end{equation}

On the other hand, from (\ref{4}) we obtain
\[
\overline{R}(e_i,e_{2m})e_{2m}=\frac{c}{4}\left[e_i+3\sum_{\beta=1}^{3}\overline{g}(e_i,J_\beta
e_{2m})J_\beta e_{2m}\right]
\]
and therefore
\begin{equation}\label{60a}
(\overline{R}(e_i,e_{2m})e_{2m})^\perp=\frac{3c}{4}\sum_{\beta=1}^{3}g(e_i,P_\beta
e_{2m})F_\beta e_{2m}.
\end{equation}

But, since $M$ is a slant submanifold, we can easily remark that
\begin{equation}\label{60b}
P_1 X=P_2 X=P_3 X,\ X\in T_pM.
\end{equation}

On the other hand, using (\ref{10b}) we obtain for all
$\beta\in\{1,2,3\}$ and $k\in\{1,...,2m\}$:
\begin{eqnarray*}
\overline{g}(F_\beta e_{2m},e_{2m+k})&=&{\rm cosec}\theta
\overline{g}(F_\beta e_{2m},F_\alpha
e_{k})\\
&=&{\rm cosec}\theta sin^2\theta g(e_{2m},e_{k})\\
&=&{\rm sin}\theta \delta_{2mk},
\end{eqnarray*}
where $\delta_{ij}$ denotes the Kronecker delta. Thus we derive
\begin{equation}\label{60c}
F_1 e_{2m}=F_2 e_{2m}=F_3 e_{2m}={\rm sin}\theta e_{4m}.
\end{equation}

From (\ref{60a}), (\ref{60b}) and (\ref{60c}) we derive
\[
(\overline{R}(e_i,e_{2m})e_{2m})^\perp=\frac{9c}{4}g(e_i,P_\alpha
e_{2m}){\rm sin}\theta e_{4m}
\]
and considering the decomposition of $P_\alpha e_{2m}$ with respect
to the adapted slant basis of $T_pM$:
\[
P_\alpha e_{2m}=\sum_{j=1}^{2m-1}\lambda_je_j
\]
we obtain
\begin{equation}\label{60}
(\overline{R}(e_i,e_{2m})e_{2m})^\perp=\frac{9c}{4}\lambda_i {\rm
sin}\theta e_{4m}.
\end{equation}

Comparing now (\ref{59}) and (\ref{60}) we derive
\[9c\lambda_i {\rm sin}\theta=0,\ i=1,...,2m-1,\]
and since $M$ is a proper slant submanifold of $\overline{M}$ and
$\sum_{i=1}^{2m-1}\lambda_i^2\neq0$, we conclude that $c=0$.
\end{proof}

\begin{cor}
There do not exist quaternionic slant immersions of minimal
codimension in $P^m(\mathbb{H})$ with unfull first normal bundle.
\end{cor}

\section*{Acknowledgements} I would like to thank the referees for carefully
reading the paper and making valuable comments
and suggestions. This work was partially supported by CNCSIS -
UEFISCSU, project PNII - IDEI code 8/2008, contract no. 525/2009.

Gabriel Eduard V\^{I}LCU$^{1,2}$ \\
      $^1$University of Bucharest,\\
      Faculty of Mathematics and Computer Science,\\
      Research Center in Geometry, Topology and Algebra,\\
      Str. Academiei, Nr. 14, Sector 1,\\
      Bucure\c sti 70109-ROMANIA\\
      e-mail: gvilcu@gta.math.unibuc.ro\\
      $^2$Petroleum-Gas University of Ploie\c sti,\\
      Department of Mathematics,\\
      Bulevardul Bucure\c sti, Nr. 39,\\
      Ploie\c sti 100680-ROMANIA\\
      e-mail: gvilcu@upg-ploiesti.ro

\end{document}